# Analytical Approximate Solutions of Systems of Multi-pantograph Delay Differential Equations Using Residual Power-series Method

Iryna Komashynska[1], Mohammed Al-Smadi [2], Abdallah Al-Habahbeh [3], Ali Ateiwi[4]

## Abstract

This paper investigates analytical approximate solutions for a system of multi pantograph delay differential equations using the residual power series method (RPSM), which obtains a Taylor expansion of the solutions and produces the exact form in terms of convergent series requires no linearization or small perturbation when the solutions are polynomials. By this method, an excellent approximate solution can be obtained with only a few iterations. In this sense, computational results of some examples are presented to demonstrate the viability, simplicity and practical usefulness of the method. In addition, the results reveal that the proposed method is very effective, straightforward, and convenient for solving a system of multi-pantograph delay differential equations.

## 1 Introduction

In this paper, we consider the following system of multi-pantograph delay differential equations

$$\begin{aligned} u_1'(t) &= \beta_1 u_1(t) + f_1\big(t, u_1(\alpha_{11}t), u_2(\alpha_{12}t), \dots, u_n(\alpha_{1n}t)\big),\\ u_2'(t) &= \beta_2 u_2(t) + f_2\big(t, u_1(\alpha_{21}t), u_2(\alpha_{22}t), \dots, u_n(\alpha_{2n}t)\big),\\ &\ \vdots\\ u_n'(t) &= \beta_n u_n(t) + f_n\big(t, u_1(\alpha_{n1}t), u_2(\alpha_{n2}t), \dots, u_n(\alpha_{nn}t)\big), \end{aligned} \tag{1}$$

subject to the initial conditions

$$u_i(t_0) = u_{i,0}\,,\ \ i = 1,2,3, \dots, n, \tag{2}$$

where $t_0 < t \leq T$, $\beta_i, u_{i,0}$ are finite constants, $f_i$ are analytical functions such that $0 < \alpha_{ij} \leq 1$, $i, j = 1,2, \dots, n$, which satisfy all necessary requirements of the existence of a unique solution, and $u_i(t)$, $i = 1,2, \dots, n$, are unknown analytical functions on the given interval to be determined.

[1]Department of Mathematics, Faculty of Science, The University of Jordan, Amman 11942, Jordan
[2]Applied Science Department, Ajloun College, Al-Balqa Applied University, Ajloun 26816, Jordan, E-Mail addresses: mhm.smadi@bau.edu.jo, mhm.smadi@yahoo.com
[3]Department of Mathematics and Computer Science, Tafila Technical University, Tafila 66110, Jordan
[4]Department of Mathematics, Faculty of Science, Al-Hussein Bin Talal University, P.O. Box 20, Ma'an-Jordan

Recently, a lot of studies about the delay differential equations (DDEs) have appeared in science literature (see e.g. [8, 22, 37]). The multi-pantograph equation is one of the most important kinds of DDEs that arise in a variety of applications in physics and engineering such as dynamical systems, electronic systems, population dynamics, quantum mechanics, etc., (for further see [18, 19, 30] and references therein). In this point, it is usually difficult to solve these kinds of DDEs analytically. Therefore, there are many powerful numerical methods in literature that can be used to approximate solutions to the multi-pantograph DDEs. To mention but a few, the Runge-Kutta method has been applied to solve multi-pantograph DDEs numerically by Li and Liu [28]. Du and Geng [20] have presented approximate solutions for singular multi-pantograph DDEs using the reproducing kernel space method (RKSM). At the same time, Yu [39] has introduced the variational iteration method (VIM) and has obtained the solution of multi-pantograph DDEs. In [36], the Taylor polynomials method has been studied for solving non-homogeneous multi-pantograph DDE with variable coefficients. In contrast, Alipour et al. [9] have used the modified VIM for finding analytical approximate solutions of multi-pantograph DDEs. In addition, Jafari and Aminataei [25] have proposed the successive approximations method for dealing with multi-pantograph DDEs and neutral functional-differential equations. Afterward, Feng [23] has employed the homotopy perturbation method (HPM) for solving multi-pantograph DDEs with variable coefficients. Lately, Geng and Qian [24] have developed a method for singularly perturbed multi-pantograph delay equations with a boundary layer at one end point based on the RKSM.

Furthermore, the RPSM has been developed as an efficient numerical as well as analytical method to determine the coefficients of power series solutions for a class of fuzzy differential equation by Abu Arqub [2]. Besides, the RPSM has been successfully applied to get numerical solutions for many other problems; For instance, generalized Lane-Emden equation which is a highly nonlinear singular differential equation [3], regular initial value problems [10] and composite and non-composite fractional differential equations [21]. This method is effective and easy to construct power series solutions for strongly linear and nonlinear equations without linearization, perturbation or discretization [6], which computes the coefficients of power series solutions by chain of linear equations of one variable. The RPSM is an alternative procedure for obtaining analytical Taylor series solution for system of multi-pantograph delay differential equations. Consequently, using the residual error concept, we get a series solution, in practice a truncated series solution. For linear problems, the exact solution can be obtained by few terms of the RPS solution. On the other hand, the numerical solvability of other version of differential problems can be found in [1, 4, 5, 7, 11-17, 27, 29, 31-35] and references therein.

The basic motivation of this paper is to apply the RPSM to develop an approach for obtaining the representation of exact and approximate solutions for system of multi-pantograph delay differential equations. This approach is simple, needs less effort to achieve the results, and effective. It does not require any converting while switching from first to higher order; thus, the method can be applied directly to the given problems by choosing an appropriate value for the initial guess approximation. Moreover, the solutions and all its derivatives are applicable for each arbitrary point in the given interval.

The remainder of this paper is organized as follows. In Section 2, basic idea of the residual power series method (RPSM) together with analysis of the method is presented. In Section 3, the RPSM is extended to provide symbolic approximate series solutions for the system of multi-pantograph equations (1) and (2). Base on the above, numerical examples are given to illustrate the capability of the proposed method in Section 3. Results reveal that only few terms are required to deduce the approximate solutions which are found to be accurate and efficient. Finally, some conclusions are summarized in the last section.

## 2 Adaptation of Residual Power Series Method (RPSM)

In this section, we present a brief description and some properties of the standard RPSM, which will be used in the remainder of this paper, in order to find out series solution for the system of multi-pantograph equations (1) and (2).

The RPSM consists of expressing the solutions of system (1) and (2) as a power series expansion about the initial point $t = t_0$. To achieve our goal, we suppose that these solutions take the following form:

$$u_i(t) = \sum_{m=0}^{\infty} u_{i,m}(t)\,, i = 1,2,\dots,n,$$

where $u_{i,m}(t)$ are the terms of approximations such that $u_{i,m}(t) = c_{i,m}(t - t_0)^m$ .

In contrast, substituting the initial guesses $u_{i,0}(\,t_0) = \frac{u_i^{(m)}(t_0)}{m!} = c_{i,0}$, which are known from initial conditions (2) for $m = 0$, into $u_i(t)$, $i = 1,2,\dots,n$, lead to the approximate solutions for the system of multi-pantograph equations $u_i(t) = u_{i,0}(\,t_0) + \sum_{m=1}^{\infty} u_{i,m}(t)$, $i = 1,2,\dots,n$, whereas $u_{i,m}(t)$, for $m = 1,2,\dots,k$, can be calculated by the following $k$th-truncated series

$$u_{i,k}(t) = \sum_{m=0}^{k} c_{i,m}(t - t_0)^m\,, i = 1,2,\dots,n. \tag{3}$$

Regarding to apply the RPSM, we rewrite the system in the following form:

$$u_i'(t) - \beta_i u_i(t) - f_i\big(t, u_1(\alpha_{i1}t), u_2(\alpha_{i2}t), \dots, u_n(\alpha_{in}t)\big) = 0, i = 1,2,\dots,n, \tag{4}$$

whereas the $k$th-residual functions and the $\infty th$ residual functions are given, respectively, by

$$\mathrm{Res}_i^k(t) = u_{i,k}'(t) - \beta_i u_{i,k}(t) - f_i\left(t, u_{1,k}(\alpha_{i1}t), u_{2,k}(\alpha_{i2}t), \dots, u_{n,k}(\alpha_{in}t)\right), i = 1,2,\dots,n, \tag{5}$$

and

$$\mathrm{Res}_i^\infty(t) = \lim_{k\to\infty} \mathrm{Res}_i^k(t) = u_i'(t) - \beta_i u_i(t) - f_i\big(t, u_1(\alpha_{i2}t), u_2(\alpha_{i2}t), \dots, u_n(\alpha_{in}t)\big), i = 1,2,\dots,n. \tag{6}$$

Obviously, $\text{Res}_i^{\infty}(t) = 0$ for each $t \in (t_0, T)$, which are infinitely differentiable functions at $t = t_0$. Furthermore, $\frac{d^m}{dt^m}\text{Res}_i^{\infty}(t_0) = \frac{d^m}{dt^m}\text{Res}_i^{k}(t_0) = 0\,, m = 0,1,2,\dots,k,$ this relation is a fundamental rule in the RPSM and its applications. Particularly, $\frac{d^{k-1}}{dt^{k-1}}\text{Res}_i^{\infty}(t_0) = \frac{d^{k-1}}{dt^{k-1}}\text{Res}_i^{k}(t_0) = 0, i = 1,2,\dots,n,\ \ k = 1,2,\dots$ .

Now, substituting the $k$th-truncated series $u_{i,k}(t)$ into Eq. (5) yields

$$\text{Res}_i^k(t) = \sum_{m=1}^{k} m c_{i,m}(t-t_0)^{m-1} - \beta_i \sum_{m=0}^{k} c_{i,m}(t-t_0)^m - f_i\left(t, \sum_{m=0}^{k} c_{1,m}(\alpha_{i1}t - t_0)^m, \sum_{m=0}^{k} c_{2,m}(\alpha_{i2}t - t_0)^m, \dots, \sum_{m=0}^{k} c_{n,m}(\alpha_{in}t - t_0)^m\right), \quad i = 1,2,\dots,n. \tag{7}$$

Consequently, based on $\text{Res}_i^1(t_0) = 0$, $i = 1,2,\dots,n$; Setting $t = t_0$, as well as $t_0 = 0$, and $k = 1$ in Eq. (7) leads to the equation

$$c_{i,1} = \beta_i c_{i,0} + f_i\big(t_0, c_{1,0}, c_{2,0}, \dots, c_{n,0}\big) = \beta_i u_{i,0} + f_i\big(0, \bar{u}_{i,0}\big), i = 1,2,\dots,n, \tag{8}$$

where $f_i\big(0, \bar{u}_{i,0}\big) = f_i\big(0, u_{1,0}, u_{2,0}, \dots, u_{n,0}\big)$. Thus, by using Eq. (3), the first approximate solutions of the system can be written as follows

$$u_{i,1}(t) = u_{i,0} + \Big(\beta_i u_{i,0} + f_i\big(0, \bar{u}_{i,0}\big)\Big)(t - t_0), i = 1,2,\dots,n. \tag{9}$$

Now, in order to obtain the second approximate solutions, we set $k = 2$ and $t_0 = 0$ such that $u_{i,2}(t) = \sum_{m=0}^{2} c_{i,m} t^m$. Then, we differentiate both sides of Eq. (7) with respect to $t$ and substitute $t = 0$ to obtain

$$\left(\frac{d}{dt}\text{Res}_i^2\right)(0) = 2c_{i,2} - \beta_i c_{i,1} - \frac{d}{dt}\left(f_i\left(t, \sum_{m=0}^{2} c_{1,m}{\alpha_{i1}}^m t^m, \sum_{m=0}^{2} c_{2,m}{\alpha_{i2}}^m t^m, \dots, \sum_{m=0}^{2} c_{n,m}{\alpha_{in}}^m t^m\right)\right), \quad i = 1,2,\dots,n.$$

According to the fact $\frac{d^{k-1}}{dt^{k-1}}\text{Res}_i^k(0) = 0, i = 1,2,\dots,n$, the values of $c_{i,2}$ are given by

$$c_{i,2} = \frac{1}{2}\left(\beta_i c_{i,1} + \frac{d}{dt}\left[f_i\left(t, \sum_{m=0}^{2} c_{1,m}{\alpha_{i1}}^m t^m, \sum_{m=0}^{2} c_{2,m}{\alpha_{i2}}^m t^m, \dots, \sum_{m=0}^{2} c_{n,m}{\alpha_{in}}^m t^m\right)\right]_{t=0}\right), i = 1,2,\dots,n.$$

Hence, the second approximate solutions of the system can be written as follows

$$u_{i,2}(t) = u_{i,0} + \left(\beta_i u_{i,0} + f_i(0, \bar{u}_{i,0})\right)(t - t_0) + \frac{1}{2}\left(\beta_i c_{i,1} + g_i(0, \bar{u}_{i,1})\right)(t - t_0)^2, i = 1,2, \dots, n, \quad (10)$$

where $g_i(0, \bar{u}_{i,1}) = \frac{d}{dt}\left[f_i\left(t, \sum_{m=0}^{2} c_{1,m}{\alpha_{i1}}^m t^m, \sum_{m=0}^{2} c_{2,m}{\alpha_{i2}}^m t^m, \dots, \sum_{m=0}^{2} c_{n,m}{\alpha_{in}}^m t^m\right)\right]_{t=0}$ and $c_{i,1}$ are given in Eq. (8) for $i = 1,2, \dots, n$.

By the same technique, the process can be repeated to generate a sequence of approximate solutions $u_{i,k}(t)$ for the system (1) and (2). Moreover, higher accuracy can be achieved by evaluating more components of the solution.

It will be convenient to have a notation for the error in the approximation $u_i(t) \approx u_{i,k}(t)$. Accordingly, we will let $\mathrm{Rem}_i^k(t)$ denote the difference between $u_i(t)$ and its $k$th Taylor polynomial; that is,

$$\mathrm{Rem}_i^k(t) = u_i(t) - u_{i,k}(t) = \sum_{m=k+1}^{\infty} \frac{1}{m!} u_i^{(m)}(t_0)(t - t_0)^m, i = 1,2, \dots, n,$$

where the functions $\mathrm{Rem}_i^k(t)$ are called the $k$th remainder for the Taylor series of $u_i(t)$. In fact, it often happens that the remainders $\mathrm{Rem}_i^k(t)$ become smaller and smaller, approaching zero, as $k$ gets large. The concept of accuracy refers to how closely a computed or measured value agrees with the truth value.

Next, we present a convergence theorem of the RPSM to capture the behavior of the solution. Afterwards, we introduce the error functions to study the accuracy and efficiency of the method. Actually, continuous approximations to the solution will be obtained.

Taylor's theorem allows us to represent fairly general functions exactly in terms of polynomials with a known, specified, and bounded error. The next theorem will guarantee convergence to the exact analytic solution of (1) and (2).

**Theorem 1:** Suppose that $u_i(x)$ is the exact solution for system (1) and (2). Then, the approximate solution obtained by the RPSM is in fact the Taylor expansion of $u_i(x)$.

***Proof.*** The proof of the Theorem is similar to proof of Theorem 1 in [2].

**Corollary 2.** Let $u_i(t), i = 1,2, \dots, n,$ be a polynomial for some $i$, then the RPSM will obtain the exact solution.

To show the accuracy of the present method for our problems, we report four types of error. The first one is called the residual error $\mathrm{Res}_i^k(t)$ and defined by

$$\mathrm{Res}_i^k(t) \coloneqq \left|u'_{i,k}(t) - \beta_i u_{i,k}(t) - f_i\left(t, u_{1,k}(\alpha_{i1}t), u_{2,k}(\alpha_{i2}t), \dots, u_{n,k}(\alpha_{in}t)\right)\right|, i = 1,2, \dots, n,$$

whilst the exact $\mathrm{Ext}_i^k(t)$, the relative error $\mathrm{Rel}_i^k(t)$ and the consecutive error $\mathrm{Con}_i^k(t)$ are defined, respectively, by

$$\mathrm{Ext}_i^k(t) := \left|u_{i,\mathrm{exact}}(t) - u_{i,k}(t)\right|,$$

$$\mathrm{Rel}_i^k(t) := \frac{\left|u_{i,\mathrm{exact}}(t) - u_{i,k}(t)\right|}{\left|u_{i,\mathrm{exact}}(t)\right|},$$

$$\mathrm{Con}_i^k(t) := \left|u_{i,k+1}(t) - u_{i,k}(t)\right|, \text{ for } i = 1,2,\dots,n,\ t \in [t_0, T],$$

where $u_{i,k}$ are the $k$th-order approximation of $u_{i,\mathrm{exact}}(t)$ obtained by the RPSM, and $u_{i,\mathrm{exact}}(t)$ are the exact solution.

## 3 Applications and Numerical Discussions

To give a clear overview of the content of this work, we consider some numerical examples to demonstrate the performance and efficiency of the RPSM. The present technique provides an analytical approximate solution in terms of an infinite power series. The consequent series truncation and the corresponding practical procedure are conducted to accomplish this task. The truncation transforms the otherwise analytical results into an exact solution that evaluated to a finite degree of accuracy. In contrast, numerical results reveal that the approximate solutions are in close agreement with the exact solutions for all values of $t$, while the accuracy is in advanced by using only few terms of approximations. Indeed, we can conclude that higher accuracy can be achieved by computing further terms. Throughout this paper, all the symbolic and numerical computations are performed using Mathematica 7.0 software package.

**Example 1**. Consider the two-dimensional pantograph equations [38]:

$$\begin{aligned} u_1'(t) - u_1(t) + u_2(t) - u_1\left(\frac{1}{2}t\right) &= f_1(t), \\ u_2'(t) + u_1(t) + u_2(t) + u_2\left(\frac{1}{2}t\right) &= f_2(t), \end{aligned} \tag{11}$$

subject to the initial conditions

$$u_1(0) = 1, u_2(0) = 1, \tag{12}$$

where $f_1(t) = e^{-t} - e^{t/2}, f_2(t) = e^{t} + e^{-t/2}$ .

To apply the RPS approach for solving system (11) and (12), we start with selecting the initial guesses of the approximations such as $u_{1,0}(t) = 1$ and $u_{2,0}(t) = 1$, then the $k$th-truncated series solutions $u_{1,k}(t)$ and $u_{2,k}(t)$ have the following form:

$$u_{1,k}(t) = \sum_{m=0}^{k} c_{1,m} t^m = 1 + c_{1,1}t + c_{1,2}t^2 + \cdots + c_{1,k}t^k,$$

$$u_{2,k}(t) = \sum_{m=0}^{k} c_{2,m} t^m = 1 + c_{2,1}t + c_{2,2}t^2 + \cdots + c_{2,k}t^k.$$

Accordingly, the unknown coefficients $c_{i,m}, m = 1,2,\dots,k, i = 1,2,$ can be found by constructing the following $k$th residual functions $\mathrm{Res}_i^k(t), i = 1,2,$ such that

$$\begin{aligned}\mathrm{Res}_1^k(t) &= \sum_{m=1}^{k} m c_{1,m} t^{m-1} - \sum_{m=0}^{k} c_{1,m} t^m + \sum_{m=0}^{k} c_{2,m} t^m - \sum_{m=0}^{k} c_{1,m}\left(\frac{t}{2}\right)^m - e^{-t} + e^{t/2},\\ \mathrm{Res}_2^k(t) &= \sum_{m=1}^{k} m c_{2,m} t^{m-1} + \sum_{m=0}^{k} c_{1,m} t^m + \sum_{m=0}^{k} c_{2,m} t^m + \sum_{m=0}^{k} c_{2,m}\left(\frac{t}{2}\right)^m - e^{t} - e^{-t/2}.\end{aligned} \tag{13}$$

Now, in order to obtain the first approximation $u_{1,1}(t)$ and $u_{2,1}(t)$ of the RPS solution for system (11) and (12), we put $k = 1$ through Eq. (13) to get

$$\mathrm{Res}_1^1(t) = c_{1,1} + \left(c_{2,1} - \tfrac{3}{2} c_{1,1}\right)\mathrm{t} - e^{-t} + e^{t/2} - 1,$$

$$\mathrm{Res}_2^1(t) = c_{2,1} + \left(c_{1,1} + \tfrac{3}{2} c_{1,1}\right)\mathrm{t} + -e^{-t/2} - e^{t} + 3\,.$$

Using the fact that $\mathrm{Res}_1^1(0) = \mathrm{Res}_2^1(0) = 0$ to get $c_{1,1} = 1$ and $c_{2,1} = -1$. Based upon this, the first approximation $u_{1,1}(t)$ and $u_{2,1}(t)$ are given by $u_{1,1}(t) = 1 + t$ and $u_{2,1}(t) = 1 - t$. Consequently, the second approximation $u_{i,2}(t), i = 1{,}2,$ of the RPS solution for system (11) and (12) can be written in the form $u_{1,2}(t) = 1 + t + c_{1,2}t^2$ and $u_{2,2}(t) = 1 - t + c_{2,2}t^2$, whereas the values of the coefficients $c_{1,2}$ and $c_{2,2}$ can be found by differentiate both sides of Eq. (13) with respect to $t$ as well as employ the RPS algorithm by putting $k = 2$ through Eq. (13) to get

$$\frac{d}{dt}(\mathrm{Res}_1^2)(t) = \left(2c_{1,2} - \tfrac{5}{2}\right) + \left(2c_{2,2} - \tfrac{5}{2} c_{1,2}\right)t + \tfrac{1}{2} e^{t/2} + e^{-t},$$

$$\frac{d}{dt}(\mathrm{Res}_2^2)(t) = \left(2c_{2,2} - \tfrac{1}{2}\right) + \left(2c_{1,2} + \tfrac{5}{2} c_{2,2}\right)t + \tfrac{1}{2} e^{-t/2} - e^{t}\,.$$

By using the fact $\frac{d}{dt}\mathrm{Res}_1^2(0) = \frac{d}{dt}\mathrm{Res}_2^2(0) = 0$, we have that $c_{1,2} = \frac{1}{2}$ and $c_{2,2} = \frac{1}{2}$. Therefore, the second approximation $u_{1,2}(t)$ and $u_{2,2}(t)$ of the RPS solution are $u_{1,2}(t) = 1 + t + \frac{1}{2}t^2$ and $u_{2,2}(t) = 1 - t + \frac{1}{2}t^2$. By continuing with the similar fashion, the 6th-order approximations $u_{i,6}(t), i = 1{,}2,$ of the RPS solution for system (11) and (12) lead to the following results:

$$\begin{aligned}u_{1,6}(t) &= 1 + t + \frac{1}{2}t^2 + \frac{1}{6}t^3 + \frac{1}{24}t^4 + \frac{1}{120}t^5 + \frac{1}{720}t^6 = \sum_{i=0}^{6} \frac{1}{(i)!} t^i,\\ u_{2,6}(t) &= 1 - t + \frac{1}{2}t^2 - \frac{1}{6}t^3 + \frac{1}{24}t^4 - \frac{1}{120}t^5 + \frac{1}{720}t^6 = \sum_{i=0}^{6} \frac{(-1)^i}{(i)!} t^i.\end{aligned} \tag{14}$$

Correspondingly, the general form of the $k$th-order RPS solutions $u_{1,k}(t)$ and $u_{2,k}(t)$ for system (11) and (12) are given by

$$u_{1,k}(t) = \sum_{m=0}^{k} c_{1,m} t^m = \sum_{j=0}^{k} \frac{1}{(j)!} t^j \text{ and } u_{2,k}(t) = \sum_{m=0}^{k} c_{3,m} t^m = \sum_{j=0}^{k} \frac{(-1)^j}{(j)!} t^j.$$

Hence, the closed forms of the RPS solutions are given by $u_1(t) = e^t$ and $u_2(t) = e^{-t}$ as soon as $k \to \infty$, which are the coinciding with the exact solutions.

To show the accuracy of the method, numerical results at some selected grid points together with comparison between the absolute errors of RPSM for 4th-order and 6th-order approximations (14) and the Laplace decomposition algorithm (LDA) [38] are given in Table 1. From the table, it can be seen that the present method provides us with an accurate approximate solution to system (11) and (12). Indeed, the results reported in this table confirm the effectiveness of the RPS method.

**Table 1:** Comparison of the absolute errors for Example 1.

| $t_i$ | Exact solution | $u_{1,k}(t_i)$ (LDA) | | $u_{1,k}(t_i)$ (Present method) | |
|---|---|---|---|---|---|
| | $u_1(t_i)$ | $k=4$ | $k=6$ | $k=4$ | $k=6$ |
| 0.2 | 1.2214027581602 | $1.210\times10^{-5}$ | $1.254\times10^{-7}$ | $2.7582\times10^{-6}$ | $2.6046\times10^{-9}$ |
| 0.4 | 1.4918246976413 | $4.238\times10^{-4}$ | $3.170\times10^{-6}$ | $9.1364\times10^{-5}$ | $3.4209\times10^{-7}$ |
| 0.6 | 1.8221188003905 | $3.499\times10^{-3}$ | $5.583\times10^{-5}$ | $7.1880\times10^{-4}$ | $6.0004\times10^{-6}$ |
| 0.8 | 2.2255409284925 | $1.594\times10^{-2}$ | $4.460\times10^{-4}$ | $3.1409\times10^{-3}$ | $4.6173\times10^{-5}$ |
| 1.0 | 2.7182818284591 | $5.236\times10^{-2}$ | $2.259\times10^{-3}$ | $9.9485\times10^{-3}$ | $2.2627\times10^{-4}$ |
| | $u_2(t_i)$ | $u_{2,k}(t_i)$ | | $u_{2,k}(t_i)$ | |
| | | $k=4$ | $k=6$ | $k=4$ | $k=6$ |
| 0.2 | 0.8187307530780 | $5.219\times10^{-5}$ | $7.807\times10^{-8}$ | $2.5803\times10^{-6}$ | $2.4776\times10^{-9}$ |
| 0.4 | 0.6703200460356 | $1.668\times10^{-3}$ | $1.310\times10^{-5}$ | $7.9954\times10^{-5}$ | $3.0952\times10^{-7}$ |
| 0.6 | 0.5488116360940 | $1.266\times10^{-2}$ | $2.227\times10^{-4}$ | $5.8836\times10^{-4}$ | $5.1639\times10^{-6}$ |
| 0.8 | 0.4493289641172 | $5.338\times10^{-2}$ | $1.668\times10^{-3}$ | $2.4044\times10^{-3}$ | $3.7791\times10^{-5}$ |
| 1.0 | 0.3678794411714 | $1.632\times10^{-1}$ | $7.956\times10^{-3}$ | $7.1206\times10^{-3}$ | $1.7611\times10^{-4}$ |

**Example 2**. Consider the system of multi-pantograph equations [38]:

$$\begin{aligned} u_1'(t) &= -u_1(t) - e^{-t}\cos\left(\frac{t}{2}\right) u_2\left(\frac{t}{2}\right) - 2e^{-(3/4)t}\cos\left(\frac{t}{2}\right)\sin\left(\frac{t}{4}\right) u_1\left(\frac{t}{4}\right), \\ u_2'(t) &= e^t u_1^2\left(\frac{t}{2}\right) - u_2^2\left(\frac{t}{2}\right), \end{aligned} \tag{15}$$

subject to the initial conditions

$$u_1(0) = 1, u_2(0) = 0\,. \tag{16}$$

Let us start with an initial approximation:

$$u_{1,0}(t) = 1,\ u_{2,0}(t) = 0. \tag{17}$$

The $k$th-truncated series formula (3) for this example by using Eq. (17) is

$$u_{1,k}(t) = \sum_{m=0}^{k} c_{1,m} t^m = 1 + c_{1,1} t + c_{1,2} t^2 + \cdots + c_{1,k} t^k,$$
$$u_{2,k}(t) = \sum_{m=1}^{k} c_{2,m} t^m = c_{2,1} t + c_{2,2} t^2 + \cdots + c_{2,k} t^k,$$

whereas the $k$th residual function $\mathrm{Res}_i^k(t), i = 1,2$, is

$$\mathrm{Res}_1^k(t) = \sum_{m=1}^{k} m c_{1,m} t^{m-1} + \sum_{m=0}^{k} c_{1,m} t^m + e^{-t} \cos\left(\frac{t}{2}\right)\left(\sum_{m=1}^{k} c_{2,m}\left(\frac{t}{2}\right)^m\right) + 2e^{-(3/4)t} \cos\left(\frac{t}{2}\right) \sin\left(\frac{t}{4}\right)\left(\sum_{m=0}^{k} c_{1,m}\left(\frac{t}{4}\right)^m\right), \tag{18}$$

$$\mathrm{Res}_2^k(t) = \sum_{m=1}^{k} m c_{2,m} t^{m-1} - e^t \left(\sum_{m=0}^{k} c_{1,m}\left(\frac{t}{2}\right)^m\right)^2 + \left(\sum_{m=1}^{k} c_{2,m}\left(\frac{t}{2}\right)^m\right)^2.$$

According to residual functions (18), the first terms of approximation of the RPS solution for $k = 1$ are $u_{1,1}(t) = 1 - t$ and $u_{2,1}(t) = t$. In contrast, the second approximation of RPS solution for this example has the form $u_{1,2}(t) = 1 - t + c_{1,2} t^2$ and $u_{2,2}(t) = t + c_{2,2} t^2$, where the values of the coefficients $c_{1,2}$ and $c_{2,2}$ can be found by differentiate both sides of Eq. (18) to construct the second-order residual functions such that

$$\begin{aligned}\frac{d}{dt}(\mathrm{Res}_1^2)(t) = {} & -1 + 2c_{1,2} + 2tc_{1,2} + 2e^{-3t/4} \cos\left(\frac{t}{2}\right) \sin\left(\frac{t}{4}\right)\left(\frac{1}{8} t c_{1,2} - \frac{1}{4}\right) \\ & + \frac{1}{2} e^{-3t/4} \cos\left(\frac{t}{4}\right) \cos\left(\frac{t}{2}\right)\left(1 - \frac{t}{4} + \frac{1}{16} t^2 c_{1,2}\right) - \frac{3}{2} e^{-3t/4} \cos\left(\frac{t}{2}\right) \sin\left(\frac{t}{4}\right) \\ & \times \left(1 - \frac{t}{4} + \frac{1}{16} t^2 c_{1,2}\right) - e^{-3t/4} \sin\left(\frac{t}{4}\right) \sin\left(\frac{t}{2}\right)\left(1 - \frac{t}{4} + \frac{1}{16} t^2 c_{1,2}\right) \\ & + e^{-t} \cos\left(\frac{t}{2}\right)\left(\frac{1}{2} + \frac{1}{2} t c_{2,2}\right) - e^{-t} \cos\left(\frac{t}{2}\right)\left(\frac{t}{2} + \frac{1}{4} t^2 c_{2,2}\right) - \frac{1}{2} e^{-t} \sin\left(\frac{t}{2}\right) \\ & \times \left(\frac{t}{2} + \frac{1}{4} t^2 c_{2,2}\right),\end{aligned}$$

$$\frac{d}{dt}(\mathrm{Res}_2^2)(t) = \frac{1}{16}\Big(-4e^t (tc_{1,2} - 1)(4 - 2t + t^2 c_{1,2}) - e^t (4 - 2t + t^2 c_{1,2})^2 + 32c_{2,2} + 4t(2 + tc_{2,2})(4 + tc_{2,2})\Big).$$

Consequently, by using $\frac{d}{dt}(\mathrm{Res}_1^2)(0) = 0, \frac{d}{dt}\left(\mathrm{Res}_2^2(0)\right) = 0$, we obtain that $u_{1,2}(t) = 1 - t$ and $u_{2,2}(t) = t$ as soon as $c_{1,2} = 0$, $c_{2,2} = 0$, which is the first approximation. Similarly, by differentiate both sides of Eq. (18) twice with respect to $t$ and using $\frac{d^2}{dt^2}(\mathrm{Res}_1^3)(0) =$ $\frac{d^2}{dt^2}\left(\text{Res}_2^2(0)\right) = 0$, the next approximation of RPS solutions is $u_{1,3}(t) = 1 - t + \frac{1}{3}t^3$, $u_2(t) = t - \frac{1}{6}t^3$ as soon as $c_{1,3} = \frac{1}{3}$ and $c_{2,3} = -\frac{1}{6}$. Furthermore, based upon $\frac{d^{k-1}}{dt^{k-1}}\ \text{Res}_i^k(0) = 0, i = 1,2,$ $k = 4,5,\ldots,10$, the 10th truncated series $u_{i,10}(t), i = 1,2$, of the RPS solution for system (15) and (16) are given as follows:

$$u_{1,10}(t) = 1 - t + \frac{t^3}{3} - \frac{t^4}{6} + \frac{t^5}{30} - \frac{t^7}{630} + \frac{t^8}{2520} - \frac{t^9}{22680},$$

$$u_{2,10}(t) = t - \frac{t^3}{6} + \frac{t^5}{120} - \frac{t^7}{5040} + \frac{t^9}{362880}.$$

Therefore, the approximate solutions of system (15) and (16) can be expressed as

$$u_1(t) = \lim_{k\to\infty} u_{1,k}(t) = \sum_{m=0}^{\infty} c_{1,m}t^m = 1 - t + \frac{t^3}{3} - \frac{t^4}{6} + \frac{t^5}{30} - \frac{t^7}{630} + \frac{t^8}{2520} - \frac{t^9}{22680} + \cdots,$$

$$u_2(t) = \lim_{k\to\infty} u_{2,k}(t) = \sum_{m=1}^{\infty} c_{2,m}t^m = t - \frac{t^3}{6} + \frac{t^5}{120} - \frac{t^7}{5040} + \frac{t^9}{362880} + \cdots,$$

that coincides with the exact solutions $u_1(t) = e^{-t}\cos t$ and $u_2(t) = \sin t$.

To illustrate the convergence of the approximate solutions $u_{i,k}(t)$ to the exact solutions $u_i(t), i = 1,2$, with respect to the $k$th-order of the solutions, we present numerical results of this example graphically. Figures 1 and 2 show the exact solution $u_i(t), i = 1,2$, and some iterated approximations $u_{i,k}(t), i = 1,2, k = 4,8,12,16,20$, respectively. These graphs reveal that the proposed method is an effective and convenient method for solving such systems with less computational and iteration steps. Moreover, in Table 2, we present numerical results with step size of 0.2 together with comparison between the absolute errors of $k$th-order, $k = 2,3$, RPS approximate solutions and LDA [38]. As a result, it is clear from this table that the approximate solutions are found to be in good agreement with the exact solutions for all values of $t$ in $[0,1]$.

**Table 2**: Comparison of the absolute errors for Example 2.

| $t_i$ | Exact solution $u_1(t_i)$ | $u_{1,k}(t_i)$ (LDA) | | $u_{1,k}(t_i)$ (Present method) | |
|---|---|---|---|---|---|
| | | $k = 2$ | $k = 3$ | $k = 2$ | $k = 3$ |
| 0.2 | 0.8024106473425 | $4.432\times10^{-4}$ | $1.900\times10^{-5}$ | $2.4107\times10^{-3}$ | $1.0647\times10^{-5}$ |
| 0.4 | 0.6174056479016 | $4.274\times10^{-3}$ | $3.656\times10^{-4}$ | $1.7406\times10^{-2}$ | $3.3898\times10^{-4}$ |
| 0.6 | 0.4529537891452 | $1.643\times10^{-2}$ | $2.119\times10^{-3}$ | $5.2953\times10^{-2}$ | $2.5538\times10^{-3}$ |
| 0.8 | 0.3130505040045 | $4.274\times10^{-2}$ | $7.420\times10^{-3}$ | $1.1305\times10^{-1}$ | $1.0650\times10^{-2}$ |
| 1.0 | 0.1987661103464 | $8.925\times10^{-2}$ | $1.960\times10^{-2}$ | $1.9877\times10^{-1}$ | $3.2099\times10^{-2}$ |
| $t_i$ | $u_2(t_i)$ | $u_{2,k}(t_i)$ (LDA) | | $u_{2,k}(t_i)$ (Present method) | |
| | | $k = 2$ | $k = 3$ | $k = 2$ | $k = 3$ |
| 0.2 | 0.1986693307951 | $5.174\times10^{-4}$ | $1.670\times10^{-5}$ | $1.3306\times10^{-3}$ | $2.6641\times10^{-6}$ |
| 0.4 | 0.3894183423087 | $5.840\times10^{-3}$ | $1.790\times10^{-4}$ | $1.0582\times10^{-2}$ | $8.5009\times10^{-5}$ |
| 0.6 | 0.5646424733950 | $2.630\times10^{-2}$ | $3.282\times10^{-4}$ | $3.5358\times10^{-2}$ | $6.4247\times10^{-4}$ |
| 0.8 | 0.7173560908995 | $8.022\times10^{-2}$ | $1.276\times10^{-3}$ | $8.2644\times10^{-2}$ | $2.6894\times10^{-3}$ |
| 1.0 | 0.8414709848079 | $1.965\times10^{-1}$ | $1.015\times10^{-2}$ | $1.5853\times10^{-1}$ | $8.1377\times10^{-3}$ |

**Example 3**. Consider the three-dimensional pantograph equations [38]:

$$u_1'(t) = 2u_2\left(\frac{t}{2}\right) + u_3(t) - t\cos\left(\frac{t}{2}\right),$$
$$u_2'(t) = 1 - t\sin(t) - 2u_3^2\left(\frac{t}{2}\right), \tag{19}$$
$$u_3'(t) = u_2(t) - u_1(t) - t\cos(t),$$

subject to the initial conditions

$$u_1(0) = -1, u_2(0) = 0, u_3(0) = 0. \tag{20}$$

Let us start with an initial approximation:

$$u_{1,0}(t) = -1,\ u_{2,0}(t) = u_{3,0}(t) = 0. \tag{21}$$

The $k$th-truncated series formula (3) for this example by using Eq. (21) is

$$u_{1,k}(t) = \sum_{m=0}^{k} c_{1,m}t^m = -1 + c_{1,1}t + c_{1,2}t^2 + \cdots + c_{1,k}t^k,$$
$$u_{2,k}(t) = \sum_{m=1}^{k} c_{2,m}t^m = c_{2,1}t + c_{2,2}t^2 + \cdots + c_{2,k}t^k,$$
$$u_{3,k}(t) = \sum_{m=1}^{k} c_{3,m}t^m = c_{3,1}t + c_{3,2}t^2 + \cdots + c_{3,k}t^k,$$

whereas the $k$th residual function $\mathrm{Res}_i^k(t), i = 1,2,3$, is

$$\mathrm{Res}_1^k(t) = \sum_{m=1}^{k} mc_{1,m}t^{m-1} - \sum_{m=1}^{k} 2c_{2,m}\left(\frac{t}{2}\right)^m - \sum_{m=1}^{k} c_{3,m}t^m + t\cos\left(\frac{t}{2}\right),$$
$$\mathrm{Res}_2^k(t) = \sum_{m=1}^{k} mc_{2,m}t^{m-1} + 2\left(\sum_{m=1}^{k} c_{3,m}\left(\frac{t}{2}\right)^m\right)^2 + t\sin(t) - 1, \tag{22}$$
$$\mathrm{Res}_3^k(t) = \sum_{m=1}^{k} mc_{3,m}t^{m-1} - \sum_{m=1}^{k} c_{2,m}t^m + \sum_{m=1}^{k} c_{1,m}t^m + t\cos(t),$$

According to residual functions (22), the first approximation of RPS solution for system (19) and (20) has the form $u_{1,1}(t) = -1 + c_{1,1}t$, $u_{2,1}(t) = c_{2,1}t$ and $u_{3,1}(t) = c_{3,1}t$, where the values of the coefficients $c_{1,1}$, $c_{2,1}$ and $c_{3,1}$ can be found by

$$\mathrm{Res}_1^1(t) = c_{1,1} - (2c_{2,1} + c_{3,1})t + t\cos\left(\frac{t}{2}\right),$$

$$\text{Res}_2^1(t) = c_{2,1} + \frac{1}{2}c_{3,1}t^2 + t\sin(t) - 1,$$

$$\text{Res}_3^1(t) = c_{3,1} + (c_{1,1} - c_{2,1})t + t\cos(t) - 1.$$

Hence, the first approximations of the RPS solution at $k = 1$ are $u_{1,1}(t) = -1$, $u_{2,1}(t) = t$ and $u_{3,1}(t) = t$. In contrast, the second approximation of RPS solution for this example has the form $u_{1,2}(t) = -1 + c_{1,2}t^2$, $u_{2,2}(t) = t + c_{2,2}t^2$ and $u_{3,2}(t) = t + c_{3,2}t^2$, where $c_{1,2}, c_{2,2}$ and $c_{3,2}$ can be found by differentiate both sides of Eq. (22) with respect to $t$ such that

$$\frac{d}{dt}(\text{Res}_1^2)(t) = 2c_{1,2} - (c_{2,2} + 2c_{3,2})t - 2 + \cos\left(\frac{t}{2}\right) - \frac{1}{2}t\sin\left(\frac{t}{2}\right),$$

$$\frac{d}{dt}(\text{Res}_2^2)(t) = 2c_{2,2} + \frac{1}{2}t(1 + tc_{3,2})(2 + tc_{3,2}) + t\cos t + \sin t + 2c_{2,2},$$

$$\frac{d}{dt}(\text{Res}_3^2)(t) = 2c_{3,2} - 2(c_{2,2} - c_{1,2})t + \cos t - t\sin t - 1.$$

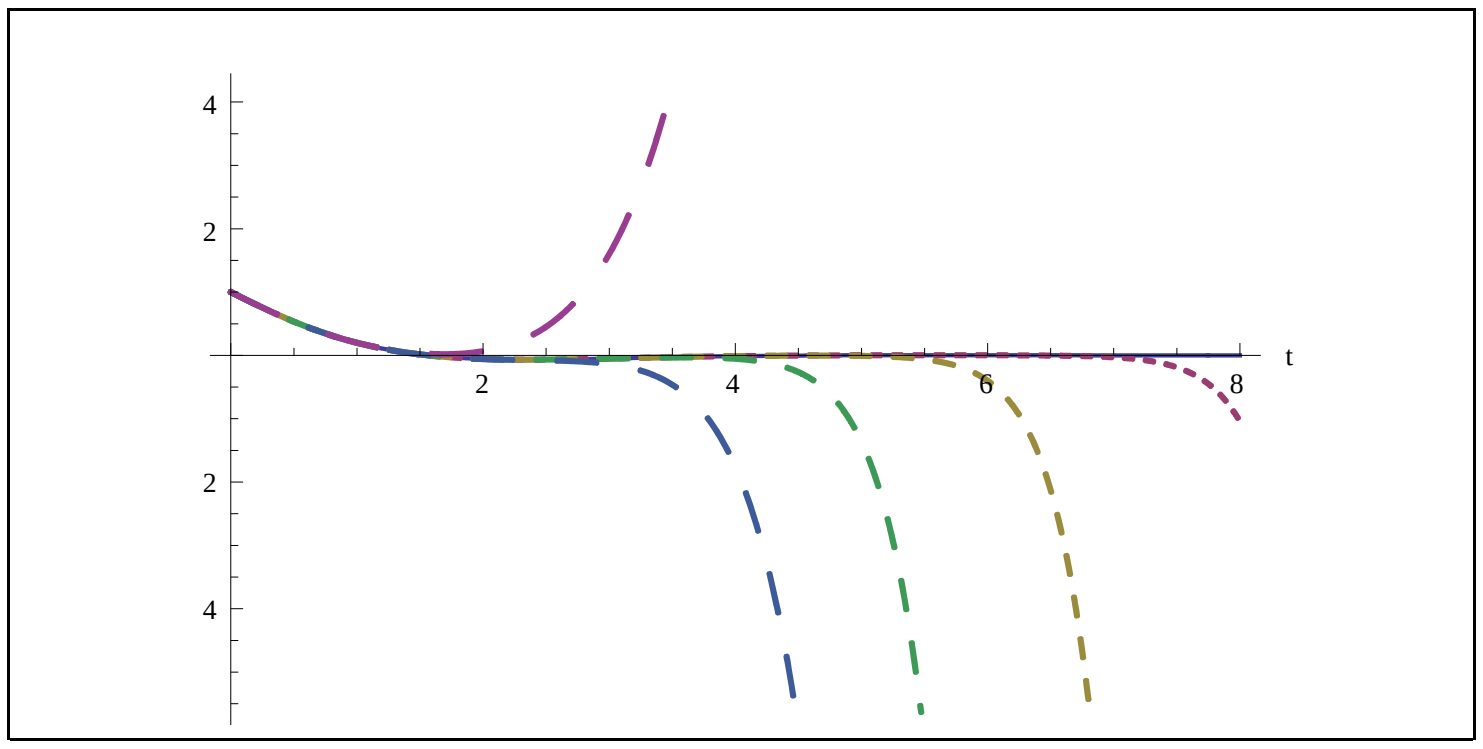

**Figure 1:** Plots of the RPS solutions $u_{1,k}(t), k = 5,10,15,20,25$, and the exact solution $u_1(t)$ of Example on $[0,8]$, where $u_1(t)$ and $u_{1,k}(t)$ are represented, respectively, by straight and dashed lines.

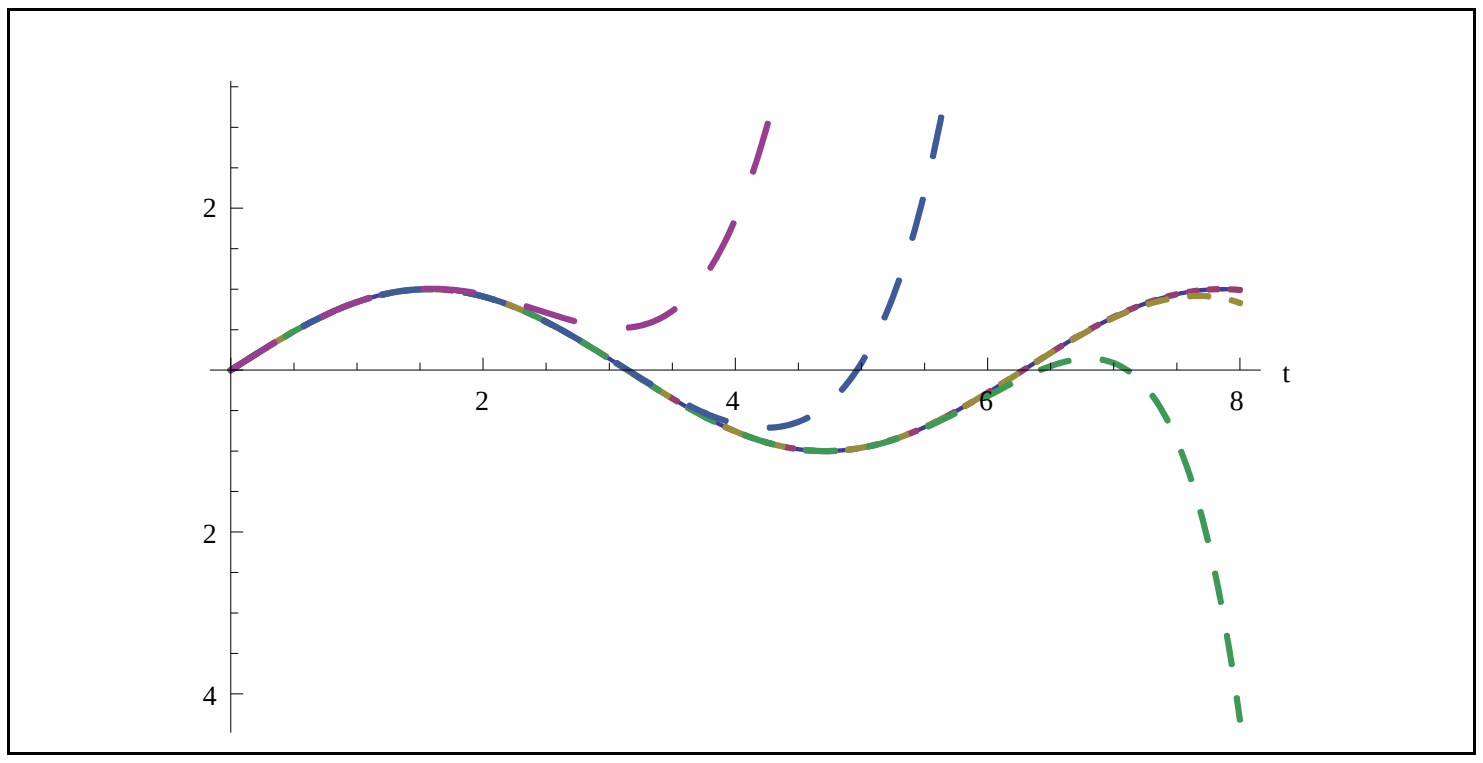

**Figure 2:** Plots of the RPS solutions $u_{2,k}(t), k = 5,10,15,20,25$, and the exact solution $u_2(t)$ of Example on $[0,8]$, where $u_2(t)$ and $u_{2,k}(t)$ are represented, respectively, by straight and dashed lines.

Thus, it is easily to get that $c_{1,2} = \frac{1}{2}$ and $c_{2,2} = c_{3,2} = 0$. Consequently, the RPS solutions of $u_{1,k}(t)$, $u_{2,k}(t)$ and $u_{3,k}(t)$ are given by

$$u_{1,k}(t) = \sum_{m=0}^{k} c_{1,m} t^m = -1 + \frac{t^2}{2} - \frac{t^4}{24} + \frac{t^6}{720} - \frac{t^8}{40320} + \cdots + c_{1,k} t^k = \sum_{j=0}^{k} (-1)^{j+1} \frac{t^{2j}}{(2j)!},$$

$$u_{2,k}(t) = \sum_{m=1}^{k} c_{2,m} t^m = t - \frac{t^3}{2} + \frac{t^5}{24} - \frac{t^7}{720} + \frac{t^9}{40320} - \cdots + c_{2,k} t^k = \sum_{j=0}^{k} (-1)^{j} \frac{t^{2j+1}}{(2j)!},$$

$$u_{3,k}(t) = \sum_{m=1}^{k} c_{3,m} t^m = t - \frac{t^3}{6} + \frac{t^5}{120} - \frac{t^7}{5040} + \frac{t^9}{362880} - \cdots + c_{3,k} t^k = \sum_{j=0}^{k} (-1)^{j} \frac{t^{2j+1}}{(2j+1)!}.$$

Hence, the closed forms of the approximate solutions as $k \to \infty$ are $u_1(t) = -\cos t$ , $u_2(t) = t\cos t$ and $u_3(t) = \sin t$ which coincides with the exact solutions.

Without loss of generality, we will test the accuracy of the present method for Example 3. Error analysis of $u_i(t), i = 1,2,3, t \in [0,1]$ for system (19) and (20) with step size of 0.2 as well as comparison among the absolute errors, relative errors, consecutive errors and residual errors of 10th-order approximate RPS solutions are shown in Tables 3, 4 and 5, respectively. From the results, it can be seen that the RPSM provides us with the accurate approximate solutions of system (19) and (20). Moreover, we can control the error also by evaluating more components of the solution.

**Table 3**: Error analysis of $u_1(t)$ for Example 3 on [0,1].

| $t$ | Absolute error $\mathrm{Ext}_1^{10}(t)$ | Relative error $\mathrm{Rel}_1^{10}(t)$ | Consecutive error $\mathrm{Con}_1^{10}(t)$ | Residual error $\mathrm{Res}_1^{10}(t)$ |
|---|---|---|---|---|
| 0.2 | 0.00 | 0.00 | 0.00 | 0.00 |
| 0.4 | $3.50830 \times 10^{-14}$ | $3.80898 \times 10^{-14}$ | $3.50830 \times 10^{-14}$ | $1.13243 \times 10^{-14}$ |
| 0.6 | $4.53548 \times 10^{-12}$ | $5.49532 \times 10^{-12}$ | $4.54448 \times 10^{-12}$ | $9.75664 \times 10^{-13}$ |
| 0.8 | $1.42961 \times 10^{-10}$ | $2.05195 \times 10^{-10}$ | $1.43464 \times 10^{-10}$ | $2.30886 \times 10^{-11}$ |
| 1.0 | $2.07625 \times 10^{-9}$ | $3.84276 \times 10^{-9}$ | $2.08768 \times 10^{-9}$ | $2.68605 \times 10^{-10}$ |

**Table 4**: Error analysis of $u_2(t)$ for Example 3 on [0,1].

| $t$ | Absolute error $\mathrm{Ext}_2^{10}(t)$ | Relative error $\mathrm{Rel}_2^{10}(t)$ | Consecutive error $\mathrm{Con}_2^{10}(t)$ | Residual error $\mathrm{Res}_2^{10}(t)$ |
|---|---|---|---|---|
| 0.2 | $5.66214 \times 10^{-15}$ | $2.88865 \times 10^{-14}$ | $5.66214 \times 10^{-15}$ | $3.10280 \times 10^{-13}$ |
| 0.4 | $1.15443 \times 10^{-11}$ | $3.13343 \times 10^{-11}$ | $1.15584 \times 10^{-11}$ | $3.17401 \times 10^{-10}$ |
| 0.6 | $9.97050 \times 10^{-10}$ | $2.01342 \times 10^{-9}$ | $9.99771 \times 10^{-10}$ | $1.82703 \times 10^{-8}$ |
| 0.8 | $2.35572 \times 10^{-8}$ | $4.22653 \times 10^{-8}$ | $2.36716 \times 10^{-8}$ | $3.23628 \times 10^{-7}$ |
| 1.0 | $2.73497 \times 10^{-7}$ | $5.06192 \times 10^{-7}$ | $2.75573 \times 10^{-7}$ | $3.00436 \times 10^{-6}$ |

**Table 5**: Error analysis of $u_3(t)$ for Example 3 on [0,1].

| $t$ | Absolute error $\text{Ext}_3^{10}(t)$ | Relative error $\text{Rel}_3^{10}(t)$ | Consecutive error $\text{Con}_3^{10}(t)$ | Residual error $\text{Res}_3^{10}(t)$ |
|---|---|---|---|---|
| 0.2 | $5.55112 \times 10^{-16}$ | $2.79415 \times 10^{-15}$ | $5.27356 \times 10^{-16}$ | $2.25653 \times 10^{-14}$ |
| 0.4 | $1.04966 \times 10^{-12}$ | $2.69546 \times 10^{-12}$ | $1.05077 \times 10^{-12}$ | $1.73516 \times 10^{-11}$ |
| 0.6 | $9.06788 \times 10^{-11}$ | $1.60595 \times 10^{-10}$ | $9.08883 \times 10^{-11}$ | $6.69236 \times 10^{-10}$ |
| 0.8 | $2.14316 \times 10^{-9}$ | $2.98758 \times 10^{-9}$ | $2.15196 \times 10^{-9}$ | $6.03226 \times 10^{-9}$ |
| 1.0 | $2.48923 \times 10^{-8}$ | $2.95819 \times 10^{-8}$ | $2.50521 \times 10^{-8}$ | $2.07625 \times 10^{-9}$ |

## 4 Conclusion

The goal of the present work was to develop an efficient and accurate algorithm to solve the system of multi-pantograph equations. This proposed algorithm produced a rapidly convergent series without need to any perturbations or other restrictive assumptions which may change the structure of the problem being solved, and with easily computable components using symbolic computation software. So, the RPSM is powerful and efficient technique in finding approximate solutions for linear and nonlinear IVPs of different types. The results obtained by the RPSM are very effective and convenient in linear and nonlinear cases because they require less computational work and time. This convenient feature confirms our belief that the efficiency of our technique will give it much greater applicability in the future for general classes of linear and nonlinear problems.

## REFERNCES